\documentclass{proc-l}

\usepackage{amssymb}

\swapnumbers
\theoremstyle{plain}
\newtheorem{theorem}[equation]{Theorem}
\newtheorem{lemma}[equation]{Lemma}
\newtheorem{corollary}[equation]{Corollary}
\newtheorem{proposition}[equation]{Proposition}

\theoremstyle{definition}
\newtheorem{definition}[equation]{Definition}
\newtheorem{example}[equation]{Example}

\numberwithin{equation}{section}

\theoremstyle{remark}
\newtheorem{remark}[equation]{Remark}

\newcommand{\complexs}{\mathbb{C}}

\newcommand{\naturals}{\mathbb{N}}

\newcommand{\integers}{\mathbb{Z}}

\newcommand{\rationals}{\mathbb{Q}}

\newcommand{\generate}[1]{\langle#1\rangle}
\newcommand{\tensor}{\otimes}
\newcommand{\into}{\hookrightarrow}

\newcommand{\iso}{\cong}

\DeclareMathOperator{\Hom}{Hom}    %Homomorphisms
\newcommand{\forget}[1]{}

\allowdisplaybreaks[2]

\newcommand{\ideal}{\mathfrak}

\newcommand{\graded}[1]{\mathrm{gr}(#1)}

\begin{document}

\title{Galois cohomology of completed link groups}
\author{Inga Blomer}
\address{Mathematisches Institut \\
Bunsenstr.~3 \\
D-37073 G\"ottingen \\
Germany}
\email{ingablomer@gmx.de}

\author[Peter Linnell]{Peter A. Linnell}
\address{Department of Mathematics \\
Virginia Tech \\
Blacksburg \\
VA 24061-0123 \\
USA}
\email{linnell@math.vt.edu}
\urladdr{http://www.math.vt.edu/people/plinnell/}

\author{Thomas Schick}
\address{Mathematisches Institut \\
Bunsenstr.~3 \\
D-37073 G\"ottingen \\
Germany}
\thanks{The third author was funded by the DAAD (German Academic
Exchange Agency)}
\email{schick@uni-math.gwdg.de}
\urladdr{http://www.uni-math.gwdg.de/schick/}

\begin{abstract}
In this paper we compute the Galois cohomology of the pro-$p$
completion of primitive link groups. Here, a primitive link group
is the fundamental group of a tame link in $S^3$ whose linking number
diagram is irreducible modulo $p$ (e.g.~none of the linking numbers is
divisible by $p$).

The result is that (with $\integers/p\integers$-coefficients) the Galois
cohomology is naturally isomorphic to the $\integers/p\integers$-cohomology of
the discrete link group.

The main application of this result is that for such groups the Baum-Connes
conjecture or the Atiyah conjecture are true for every finite extension (or
even every elementary amenable extension), if
they are true for the group itself.
\end{abstract}

\subjclass[2000]{Primary: 20E18; Secondary: 20J06, 57M25}

\keywords{Link group, Lie algebra, Galois cohomology}

\maketitle

\section{Introduction}

Throughout the paper we fix a prime number $p$. If $G$ is a discrete
group, then $\hat G$ will denote its pro-$p$ completion.
  
 \begin{definition} We denote the fundamental group $G$ of the
  complement $M$ of a tubular neighborhood of a tame link $L$ with $d$
  components in
  $S^3$  a \emph{link group with $d$ components}.
  We define the  \emph{linking diagram} to be the edge-labeled graph whose
  vertices
  are the components of the link. Any two vertices are joined by one
  edge, which is labeled with the linking number of the two link
  components involved.

  We say the link group $G$ is \emph{$q$-primitive} for a prime $q$ if
  there is a spanning subtree of the linking diagram such that none of
  the labels of the edges of this subtree is congruent to $0$ modulo
  $p$. It is called \emph{primitive} if it is $q$-primitive for every
  prime $q$.

  Observe that in particular every knot group (i.e.~a link group of a
  link with only one component) is primitive, because there are no edges to
  worry about.
\end{definition}

Labute computes in \cite{Labute(1990)} the lower central series
quotients of primitive link groups. He shows in particular that all of
them are torsion free. He uses Lie algebra techniques
developed in \cite{Labute(1967),Labute(1970),Labute(1985)}. These
techniques give results not only about the lower central series quotients, but
also about the pro-$p$ completions and their Galois cohomology.

In \cite{Linnell-Schick(2000b),Schick} it became important to compute the
Galois cohomology of such groups because of the following result
which combines special cases of \cite[Theorem 4.60, Corollary
4.62]{Linnell-Schick(2000b)} or of \cite[Theorem 10]{Schick}. 

\begin{theorem}\label{theorem:generalAtiyah}
  Let $G$ be a discrete group with finite classifying space. Assume that all the lower central
  series quotients of $G$ are torsion free and that for every prime
  $q$ the canonical map
  \begin{equation}\label{eq:1a}
    H^*({\hat G},\integers/q) \to H^*(G,\integers/q)
  \end{equation}
  between the Galois cohomology of the pro-$q$ completion $\hat G$ and the
  ordinary cohomology of $G$ is an isomorphism.

  If the Atiyah conjecture holds for $G$, then it holds also for
  every elementary amenable extension of $G$, and in particular every
finite extension of $G$.

  If the Baum-Connes conjecture (with coefficients) holds for $G$, then it
  holds also for every elementary amenable extension of $G$, and in
particular every finite extension of $G$.
\end{theorem}

Here, the Atiyah conjecture for a torsion-free group $G$ asserts that
the $L^2$-Betti numbers of every compact manifold with fundamental
group $G$ are integers. This implies that the ring $\rationals G$
embeds into a skew field (there are versions which imply the same for
$\complexs G$).

It is checked in the proof of \cite[Proposition
5.34]{Linnell-Schick(2000b)} that a primitive link group fulfills all the
assumptions of Theorem 
\ref{theorem:generalAtiyah}, except for the computation
of the Galois cohomology. This will be done in this paper. Our main
result is the following.
\begin{theorem}\label{theorem:main}
  Let $G$ be a primitive link group. Then the map
  \begin{equation*}
    H^n({\hat G},\integers/p) \to H^n(G,\integers/p)
  \end{equation*}
  is an isomorphism for every $n\in\integers$ and every prime number
$p$.
\end{theorem}
\begin{definition}\label{def:cohom_complete}
  As in \cite[Definition 4.3]{Linnell-Schick(2000b)}, we call the above
  property, that $H^n({\hat G},\integers/p) \to H^n(G,\integers/p)$ is an
  isomorphism for every $n\in\integers$ and every prime number $p$, as $G$
  being \emph{cohomologically complete}.
\end{definition}

\begin{corollary}
Let $G$ be a primitive link group.
  If $G$ satisfies the Atiyah conjecture, then so does every
  elementary amenable extension of $G$, and in particular every
finite extension of $G$.

  The Baum-Connes conjecture is true for every elementary amenable
  extension of $G$, and in particular every finite extension of $G$.
\end{corollary}
\begin{proof}
  All these statements follow immediately from Theorems \ref{theorem:main} and
  \ref{theorem:generalAtiyah}, and from the fact that the Baum-Connes
  conjecture with coefficients is known for fundamental groups of Haken
  $3$-manifolds, compare e.g.~\cite[Theorem 5.23]{MR2027168} or
  \cite[Theorem 5.2]{MR2181833}.
\end{proof}

This work was inspired by the confirmation of John Labute that
the homomorphism of
(\ref{eq:1a}) indeed should be an isomorphism for primitive link groups, and it
is very much based on Labute's work cited above. Since the 
result seems not to be available in the literature, and since it is not
obvious how to prove it using
\cite{Labute(1967),Labute(1970),Labute(1985),Labute(1990)} we indicate
the proof in this note. We will heavily refer to results and methods
which can be found in these papers. Alternatively, one could have used a
suitable version of the mildness and appropriate versions of the results of
\cite{Labute05}. However, since there the $p$-central series is in the focus,
whereas we use the central series, we try to give a somewhat independent
treatment.

\section{Proof of Theorem \ref{theorem:main}}
\label{sec:proof}

First we observe (cf.~\cite[Lemma 4.5]{Linnell-Schick(2000b)}):
\begin{lemma}\label{lem:low_completeness}
  Let $G$ be a discrete group. Then
  \begin{equation*}
    H^k({\hat G},\integers/p) \to H^k(G,\integers/p)
  \end{equation*}
  is an isomorphism for $k=0$ and $k=1$.
\end{lemma}
\begin{proof}
  For $k=0$ this is immediately clear from the definition. Also,
  $H^1(G,\integers/p)=\Hom(G,\integers/p)$ and $H^1({\hat
  G},\integers/p)=\Hom({\hat G},\integers/p)$. But by the very
  definition of the pro-$p$
  completion every homomorphism from $G$ to $\integers/p$ extends
  uniquely to a homomorphism  ${\hat G}\to \integers/p$.
\end{proof}

\begin{example}\label{ex:free_completeness}
  Let $F$ be an arbitrary free group (not necessarily finitely
  generated). Then $F$ is cohomologically complete, as defined in Definition
  \ref{def:cohom_complete}. 
\end{example}
\begin{proof}
  Since $F$ has cohomological dimension $\le 1$,
  in view of Lemma \ref{lem:low_completeness} it suffices to prove
  that the pro-$p$ completion of $F$ has cohomological dimension $\le
  1$, i.e.~that ${\hat F}^p$ is a free pro-$p$ group. But this follows
  from \cite[Corollary 4 and Remark on p.~8]{Serre(1997)}.
\end{proof}

The first step to prove Theorem \ref{theorem:main} is the computation
of the cohomology of the discrete link group:
\begin{lemma}\label{lem:linkcohom}
  Assume that $G$ is a primitive link group with $d$ components. Then $G$
  has a compact $2$-dimensional classifying space and 
  \begin{equation*}
\dim_{\integers/p\integers}
  H^n(G,\integers/p\integers)=
  \begin{cases}
    1; & \text{if }n=0\\ d; & \text{if } n=1\\ d-1; & \text{if }n=2 \\
    0; & \text{if }n>2.
  \end{cases}
\end{equation*}
\end{lemma}
\begin{proof}
  Let $L$ be a link with $d$ components with primitive linking
  diagram such that $G$ is the fundamental group of $S^3-L$. It is
  shown in the proof of \cite[Proposition 5.34]{Linnell-Schick(2000b)}
  that the compact $3$-dimensional manifold with boundary $S^3-U(L)$
  ($U(L)$ being an open tubular neighborhood of $L$) is a classifying
  space for $G$, homotopy equivalent to a $2$-dimensional CW-complex. In
  particular,
  $H^n(G,\integers/p\integers)=H^n(S^3-L,\integers/p\integers)$ which easily
  can be
  computed using Alexander duality and gives the asserted dimensions.
\end{proof}

These computations show that the following proposition is relevant for
$p$-primitive link groups.

\begin{proposition}\label{prop:two_dim_completeness}
  Let $G$ be a discrete group with finite classifying space of
  dimension $2$. Then 
  \begin{equation}\label{eq:2a}
    H^n({\hat G},\integers/p) \to H^n(G,\integers/p)\quad\text{is
    an isomorphism for every $n\in\integers$}
  \end{equation}
  if and only if 
  \begin{enumerate}
  \item\label{item:1} $H^3({\hat G},\integers/p\integers) =0$
  \item \label{item:2} $\dim_{\integers/p\integers} H^2(\hat{G},
\integers/p\integers) =
    \dim_{\integers/p\integers} H^2(G,\integers/p\integers)$.
  \end{enumerate}
  Condition \ref{item:1} is equivalent to the fact that the
  cohomological dimension of ${\hat G}$ is $\le 2$. If ${\hat G}={
  F}/ R$ where
  $F$ is a finitely generated free pro-$p$ group, this is also
  equivalent to $
  R/\overline{[R,R]}$
  (where $\overline{[R,R]}$ is the closure of the subgroup generated by all
  commutators of $R$) being a free $\integers_p[[{\hat
  G}]]$-module. $\integers_p[[{\hat G}]]$ is the completed group
  ring (compare \cite{Brumer(1966)}) and the action is induced by
  conjugation.
\end{proposition}
\begin{proof}
  \eqref{eq:2a} obviously implies that \ref{item:1}
  and \ref{item:2} hold. 

  By \cite[Propositions 11 and 20]{Serre(1997)} we have $H^3({\hat
  G},\integers/p)=0$ if and only if the cohomological dimension of
  ${\hat G}$ is $\le 2$. The equivalence of the last condition with
  \ref{item:1} follows immediately from \cite[Theorem 5.2 and Corollary
  5.3]{Brumer(1966)}. Here we use the fact that the abelian group $N_p:=
  R/\overline{[R,R]}$
  is equal to its
  $p$-Sylow subgroup since a free pro-$p$ group (like $R$)
  does not have non-trivial cyclic quotients of prime order except for
  $\integers/p$,
  so that our definition indeed coincides with the definition of
  \cite{Brumer(1966)}.

  We use \cite[Lemma 4.2]{Brumer(1966)} to compute the Galois
  cohomology of ${\hat G}$ using a $\integers_p[[{\hat
  G}]]$-resolution of the $p$-adic integers $\integers_p$.

  We change notation and write $G=F/R$ with a free group $F$ on $d$
 generators $x_1,\dots,x_d$ and
 relations $R$ (normally generated by $r$ relators $r_1,\dots, r_r$, and where
 $r$ is minimal with this property), where the
 presentation is given by a finite
 classifying $2$-complex $X$ for $G$.

 Then ${\hat 
 G} = {\hat F}/\overline R$, where $\overline R$ is the closed
 normal subgroup of ${\hat F}$ generated by $R$ (compare \cite[Lemma 3.9]{Linnell-Schick(2000b)}).

 In the
  proof of \cite[Theorem 5.2]{Brumer(1966)} the exact sequence
  \begin{equation}\label{eq:1}
    0\to N_p \to I(F)\hat\tensor_{\integers_p[[N_p]]}\integers_p \to \integers_p[[{\hat G}]] \to
    \integers_p\to 0 
  \end{equation}
  is given. Since $F$ has cohomological dimension $1$, the
  discussion at the beginning of \cite[Proof of Corollary
  5.3]{Brumer(1966)} shows that $I(F)$ is a free $\integers_p[[{\hat
  F}]]$-module, generated by $\{(1-x_1),\dots,(1-x_d)\}$ (note that $\{x_1,\dots
  x_d\}$ constitute also a set of free
generators for ${\hat F}$). This fact is also used in \cite[Example (2)
 on p.~144]{Labute(1967)}.  Hence
  $I(F)\hat\tensor_{\integers_p[[N_p]]}\integers_p
\iso \integers_p[[{\hat G}]]^d$.

As mentioned by Brumer in \cite{Brumer(1966)},
his reasoning  applies also to the discrete
group $G$ to give the exact sequence
\begin{equation}\label{eq:2}
  0 \to N \xrightarrow{d_2} \integers[G]^d \xrightarrow{d_1} \integers[G]
  \xrightarrow{d_0} \integers \to 0,
\end{equation}
where $N=R/[R,R]$ is the abelianization of the relations. It follows
from the construction that the inclusions $G\to{\hat G}$ and
$\integers\into \integers_p$ induce maps
which are just the inclusion of a free $\integers[G]$-module into its
completion for the two middle terms in the exact sequences. Moreover,
the boundary map $d_1\colon \integers[G]^d\to\integers[G]$ maps the
generator of copy number $j$ of $\integers[G]$ to $1-x_j$, where $x_j$
is (the image of) one of the generators of $F$.

On the other hand, the (universal covering of the) $2$-complex $X$ has
a cellular chain complex 
\begin{equation*}
  0\to  \integers[G]^r \xrightarrow{d_2^X} \integers[G]^d \xrightarrow{d_1^X}\integers[G]\xrightarrow{d^X_0}
  \integers\to 0.
\end{equation*}
Since this universal covering is contractible by assumption, this
sequence is exact (this is the reason why we added the group
$\integers$ at the end). Moreover, the differentials $d_1^X$ and $d_1$
as well as $d_0^X$ and $d_0$ coincide. Identifying $N$ and
$\integers[G]^r$ with their images under $d_2$ or $d_2^X$,
respectively, we see that they are isomorphic, in particular $N$
is a free $\integers[G]$-module, freely generated by the relators
$r_1,\dots, r_r$. To compare the sequence \eqref{eq:2} with
\eqref{eq:1} we observe that by \cite[Corollary
5.3]{Brumer(1966)} (under the assumption that ${\hat G}$
has cohomological dimension $\le 2$) $N_p$ is a free
$\integers_p[[{\hat G}]]$-module, but not necessarily generated by
all the $r_1,\dots, r_r$, instead only by a minimal subset which still
normally and topologically generates $\overline{R}$ (which a
priori might be smaller). By \cite[Proof of Corollary 5.3]{Brumer(1966)}, we
need any subset whose image in $\overline{R}/\overline{[R,F]R^p}$ is a basis
as $\integers/p\integers$-vector space.

To compute the cohomology with $\integers/p\integers$-coefficients of $G$ or
${\hat G}$, respectively, we have to pass to the (topological) dual
chain complexes and take their cohomology.

By the very definition of the completion,
$\Hom(\integers[G]^d,\integers/p\integers)$ is isomorphic to
$\Hom(\integers_p[[{\hat G}]]^d,\integers/p\integers)$ (every such homomorphism
extends uniquely). However, the same considerations only show that
$\Hom(N_p,\integers/p\integers)\to  
\Hom(N,\integers/p\integers)$ is injective. Diagram chasing reveals that the
induced maps $H^j({\hat G},\integers/p\integers)\to
H^j(G,\integers/p\integers)$ are 
isomorphisms for $j=0,1$ 
and injective for $j=2$. Since we assume that for $j=2$ the two
cohomology groups have the same finite dimension, the map is an
isomorphism in this case, as well.
\end{proof}

To apply Proposition \ref{prop:two_dim_completeness}, we next derive a
presentation of the pro-$p$ completion of link groups.

\begin{proposition}\label{prop:linkpresentation}
  Let $G$ be a link group with $d$ components. Then the pro-$p$
  completion ${\hat G}$ of $G$ has a presentation as a pro-$p$ group
  with $d$ generators and $d-1$ relators. In other words, ${\hat G}$
  is a quotient of a free pro-$p$ group $F$ with $d$ generators, and the
  kernel $R$ is the closure of the normal subgroup generated by $d-1$
  elements.

  If, moreover, $G$ is primitive and 
  if the $\integers_p[[{\hat G}]]$-module $R/\overline{[R,R]}$ is
  free on $d-1$ generators then the natural map
  \begin{equation*}
    H^n({\hat G},\integers/p\integers)\to H^n(G,\integers/p\integers)
  \end{equation*}
  is an isomorphism for every $n\in\integers$.
\end{proposition}
\begin{proof}
  By a result of Milnor (compare
  \cite[p.~951]{Labute(1990)}), in our situation the lower central
  series quotients $G/\gamma_n(G)$ have presentations
  \begin{equation*}
    \generate{x_1,\dots,x_d \mid
    [x_1,w_1^{(n)}]=\dots=[x_{d-1},w_{d-1}^{(n)}]=1,
    \gamma_n(F)} 
\end{equation*}
where $F$ is the free group generated by $x_1,\dots,x_d$,
$[x,y]=x^{-1}y^{-1}xy$ is the commutator, and for $k=1,\dots,d-1$ and
  $n\in\integers$ we have certain
$w_k^{(n)}\in F$  which in addition fulfill $w_k^{(n)} \equiv w_k^{(n+1)}$
modulo $\gamma_n(F)$.

Consequently, if $\gamma_n^p(G)$ is the lower
$p$-central series (defined inductively by $\gamma_1^p(G)=G$ and
$\gamma^p_{n+1}(G)=(\gamma_n(G))^p[\gamma^p_n(G),G]$) then
$w_k^{(n)}\equiv w_k^{(n+1)}$ modulo $\gamma_n^p(F)$, and
$G/\gamma_n^p(G)$ has the presentation
\begin{equation*}
    \generate{x_1,\dots,x_d \mid
    [x_1,w_1^{(n)}]=\dots=[x_{d-1},w_{d-1}^{(n)}]=1,
    \gamma_n^p(F)}.
\end{equation*}
Since ${\hat G}$ is the inverse limit of the inverse system given by
$G/\gamma_n^p(G)$  it follows
that for each $k=1,\dots,d$ the sequence $w_k^{(n)}$ converges in
${\hat G}$ to $w_k$ and that a presentation of ${\hat G}$ as
pro-$p$ group is given by
\begin{equation*}
  \generate{x_1,\dots,x_d \mid [x_1,w_1]=\dots=[x_{d-1},w_{d-1}]=1}.
\end{equation*}

Assume now that $G$ is primitive. By Lemma \ref{lem:linkcohom} $G$ has
a classifying space of dimension $2$. Because of the assumptions and
Proposition \ref{prop:two_dim_completeness} we only have to check that
$\dim_{\integers/p} H^2(G,\integers/p)= \dim_{\integers/p} H^2({\hat
  G},\integers/p)$. In the proof of Proposition
\ref{prop:two_dim_completeness} we have seen that the cohomology of
${\hat G}$ is the cohomology of the following cochain complex (we
dualize \eqref{eq:1} using the fact that $I(F)$ is a free
$\integers_p[[{\hat G}]]$-module on $d$ generators and by assumption
$N_p=R/\overline{[R,R]}$ is free on $d-1$ generators)
\begin{equation}\label{eq:3}
0 \leftarrow
  (\integers/p\integers)^{d-1}\xleftarrow{d_1^*} (\integers/p\integers)^d
  \xleftarrow{d_0^*} \integers/p\integers.
\end{equation}
By Lemma \ref{lem:low_completeness} and Lemma \ref{lem:linkcohom} we
know that $\dim_{\integers/p\integers} H^0({\hat G},\integers/p\integers)=1$ and
$\dim_{\integers/p\integers} H^1({\hat G},\integers/p\integers)=d$. This implies that
$d_1^*$ in the cochain complex \eqref{eq:3} is trivial, and therefore
$\dim_{\integers/p\integers}H^2({\hat G},\integers/p\integers)=d-1 =\dim_{\integers/p\integers}
H^2(G,\integers/p\integers)$.

Note that it is crucial not only to know that $R/\overline{[R,R]}$ is
free, but that it is free on $d-1$ free generators.
\end{proof}

\begin{remark}
  It should be observed that it is not true in general that a discrete
  link group has a presentation with $d$ generators and $d-1$
  relators, so that the presentation obtained in Lemma
  \ref{prop:linkpresentation} will in general not come from a
  presentation of $G$ itself.
\end{remark}

In order to show that primitive link groups are cohomologically complete, it
remains to show that for a primitive $d$-component link group $G$ with pro-$p$
completion 
$\hat{G}=F/R$ the following holds:  
\[R/\overline{[R,R]} \text{ is a free }\integers _p[[\hat{G}]]\text{-module on
  $(d-1)$-generators.}\]

% Labute has shown that this holds for so-called \emph{mild} pro-$p$ groups
% \cite[Theorem 
% 1.2 and Theorem 4.1]{Labute05}.

% We next recall this concept and will state precisely the main theorem of Labute
% about mild groups we are going to use, and then will show that a primitive
% link group $G$ indeed is mild in such a way that we can apply Labute's result
% to prove Theorem \ref{theorem:main}.

From now on let $F$ be a finitely generated free pro-$p$ group and $R$
a closed subgroup of $F$. Set $G:=F/R$.
We now give a criterion in terms of filtrations, graded rings and
graded modules for $N:=R/\overline{[R,R]}$ to be a free
$\integers_p[[G]]$-module.

Observe that the lower central series $\gamma_1(F)=F$,
$\gamma_k(F)=\overline{[\gamma_{k-1}(F),F]}$ provides a filtration of
$F$  by closed subgroups, which induces a filtration $N_k$ of $N$,
where $N_k$ is the image of $R\cap\gamma_k(F)$ under the projection
$R\to N$. 

This induces graded groups $\graded{F}$ with $\graded{F}_k
:=\gamma_k(F)/\gamma_{k+1}(F)$. Each $\graded{F}_k$ is an abelian
pro-$p$ group, i.e.~a $\integers_p$-module. Moreover, the commutator
induces a Lie ring structure on the direct sum $\graded{F}$ of all
the $\graded{F}_k$, so that $\graded{F}$ is a $\integers_p$-Lie
algebra (compare \cite[Section 2.1]{Labute(1966)}). By
\cite[Proposition 1]{Labute(1966)} $\graded{F}$ is a free Lie algebra
over $\integers_p$ generated by the image of the free (topological)
generators of $F$. Similarly, we define $\graded{N}$ by
$\graded{N}_k:=N_k/N_{k+1}$. Observe that $\graded{N}$ inherits a
  $\integers_p$-Lie algebra structure, but it is an abelian Lie
  algebra (i.e.~the bracket is trivial).

Let $I$ be the augmentation ideal of $\integers_p[[G]]$, i.e.~the
kernel of the augmentation homomorphism $\integers_p[[G]]\to
\integers_p$. Then the closures of the powers of $I$ provide a
filtration of $\integers_p[[G]]$ by closed ideals, and we get the
graded ring $\graded{\integers_p[[G]]}$ with
$\graded{\integers_p[[G]]}_k=\overline{I^k}/\overline{I^{k+1}}$ ($k=0,1,\dots$). It is easy to check that
the $\integers_p[[G]]$-module structure of $N$ (induced from
conjugation of $G$ on $R$) fulfills $\overline{I^n}N_k\subset
N_{n+k}$. Therefore, it induces a $\graded{\integers_p[[G]]}$-module
structure on $\graded{N}$. 

\begin{lemma}\label{lem:graded_to_ungraded_free}
  Assume, in the above situation, that $w_1,\dots,w_r$ are elements of $N$ with
  images $W_1,\dots,W_r$ in $\graded{N}$. Assume further that $\graded{N}$ is a
  free $\graded{\integers_p[[G]]}$-module with basis
  $W_1,\dots,W_r$. Then $N$ is a free $\integers_p[[G]]$-module with
  basis $w_1,\dots,w_r$.
\end{lemma}
\begin{proof}
  The elements $w_1,\dots,w_r$ determine a $\integers_p[[G]]$-module map
  $\phi\colon \integers_p[[G]]^r\to N$. We want to show that $\phi$ is an
  isomorphism.

  An easy induction and the assumptions imply that the induced map
  \begin{equation*}
\phi_k\colon (\integers_p[[G]]/\overline{I^{k-1}})^r \to N/N_k
\end{equation*}
is an
  isomorphism for every $k\ge 1$. It follows from \cite[Lemme
  1]{Douady(1960)} that $\phi$ induces an isomorphism
  \begin{equation}\label{eq:4}
    \phi_\infty\colon (\integers_p[[G]]/I^\infty)^r \xrightarrow{\iso} N/N_\infty
  \end{equation}
  where $I^\infty=\bigcap_{k\in\naturals} \overline{I^k}$ and
  $N_{\infty}=\bigcap_{k\in\naturals} N_k$ (here $\naturals$ denotes
the positive integers $\{1,2,\dots\}$).

  Let $J$ be the radical of $\integers_p[[G]]$ (i.e.~the kernel of the
  projection $\integers_p[[G]]\to\integers/p\integers$). Then certainly
  $I\subset J$. By \cite[Chap.~II (2.2.2)]{Lazard(1965)} $J^\infty$ is
  trivial, therefore the same is true for $I^\infty$.

  It is well known that $\bigcap_{k\in\naturals}\gamma_k(F)=\{1\}$ and
  therefore $\bigcap_{k\in\naturals} (R\cap\gamma_k(F))
  =\{1\}$. Assume that $w\in N_\infty$. Since
  $R\cap\gamma_k(F)$ maps onto $N_k$, for each $k\ge 1$ we find
  $r_k\in R\cap \gamma_k(F)$ such that the image of $r_k$ in $N_k$ is
  equal to $w$. Compactness and completeness of $F$ (and $R$) imply
  that we find a subsequence of $r_k$ which converges to $r\in
  R$. Since $r_j\in\gamma_k(F)$ for $j\ge k$ and since $\gamma_k(F)$
  is closed, $r\in\gamma_k(F)$ for every $k$. It follows that $r$ is
  trivial. Because of continuity of the projection map $R\to N$, the
  image of $r$ in $N$ is $w$. Hence $w$ and therefore $N_\infty$ is
  trivial.

  Therefore \eqref{eq:4} proves the lemma.
\end{proof}

The definition of the lower central series of a group and the
corresponding graded ring makes sense also in the discrete case, and
similarly with all the other objects defined above. We get a parallel
theory, with the completed group ring replaced by the group ring over
$\integers$. The discrete situation was studied thoroughly
by Labute. We now give a connection between the discrete and the
pro-$p$ completed versions. In particular, this will give a condition
when $\graded{N}$ is a free $\graded{\integers_p[[G]]}$-module.

To do this, we recall the following lemma, whose proof we give for the
convenience of the reader.

\begin{lemma}\label{lem:unique_completion}
  Let $\Gamma$ be a discrete group and $i\colon\Gamma\to G$ be an injective
  group homomorphism of $\Gamma$ into a pro-$p$ group $G$ with dense
  image. Then $G$ is the pro-$p$ completion of $\Gamma$.
\end{lemma}
\begin{proof}
  By the universal property of the pro-$p$ completion ${\hat\Gamma}$
  we have a factorization of $i$ as
  $\Gamma\to{\hat\Gamma}\xrightarrow{I} G$. Here $I$ is surjective,
  since the image is closed and dense. We identify $\Gamma$ with its
  image in ${\hat\Gamma}$ and in $G$. It remains to check
  injectivity of $I$. Note that every homomorphism of $\Gamma$ onto a finite
  $p$-group extends uniquely to a homomorphism of ${\hat\Gamma}$. It
  follows in particular that for a closed normal subgroup $U$ of
  $\hat\Gamma$ of finite index $p^n$ the intersection with $\Gamma$
  has index $p^n$ in $\Gamma$. Since $I$ is surjective and $U$ is
  compact, $I(U)$ is a closed normal subgroup of $G$ of index not
  bigger than $p^n$. Since $i$ is injective, $I(U)\cap\Gamma$ has
  index $p^n$ in $\Gamma$. It follows that the index of $I(U)$ in $G$
  is precisely $p^n$. Therefore the induced map $I_U:\hat\Gamma/U\to
  G/I(U)$ is an isomorphism. By \cite[Lemme 1]{Douady(1960)} $I$
  induces an isomorphism of $\hat\Gamma$
  onto a quotient of $G$. It follows that $I$ itself is an isomorphism.
\end{proof}

  Let $L$ be a free $\integers$-Lie algebra on $d$ generators. It is a
  classical fact that this
  is the graded Lie algebra of the lower central series of a discrete free
  group $\Gamma_d$ on $d$ 
  generators (compare \cite[Theorem 6.1]{Serre(1965)}). If
  $F$ is a free pro-$p$ group on $d$ generators then we have a
  canonical map $L\to \graded{F}$ (induced from the map $\Gamma_d\to
  F$ given by the generators) and since $\graded{F}$ is a free
  $\integers_p$-Lie algebra on the given set of generators, for each
  $k\in\naturals$ 
  $L_k\to \graded{F}_k$ is just the pro-$p$ completion of the
  (finitely generated, compare e.g.~\cite[Corollary
  5.12]{Magnus-Karras-Solitar(1966)}) abelian
  group $L_k$. It is also possible to conclude this from the fact
  that $\Gamma_d/\gamma_k(\Gamma_d)$ injects into its pro-$p$
  completion by \cite[Theorem 2.2]{Gruenberg(1957)} together with Lemma
  \ref{lem:unique_completion}. 

The bracket $L_k\times L_r\to L_{k+r}$
  extends by 
  continuity to $\graded{F}_k\times \graded{F}_r\to
  \graded{F}_{k+r}$. This follows because the bracket is induced from
  taking commutators, i.e.~from group multiplication, and the latter is
  continuous (we then only have to pass to quotients). 

\begin{lemma}\label{lem:ideal_closure}
  Let $w_1,\dots,w_r$ be elements of $F$ which topologically generate a closed
  normal subgroup $R$ of $F$. Let $W_1,\dots,W_r$ be the natural
  images in $\graded{F}$. Assume that the set $\{W_1,\dots,W_r\}$ belongs to
  the 
   $\integers$-Lie subalgebra $L$ of $\graded{F}$. Let $\ideal{r}$ be the
  (graded) $\integers$-Lie ideal 
  of $L$ generated by $(W_1,\dots,W_r)$, and let $\ideal{R}$ be the
  closed (graded) $\integers_p$-Lie ideal of $\graded{F}$ generated by
  $(W_1,\dots,W_r)$
(closed means each $\ideal{R}_k$ is closed in $\graded{F}_k$). 
 
  In this situation, $\ideal{R}_k$ is the closure of $\ideal{r}_k$
  in $\graded{F}_k$.
  Moreover, $(\graded{F}/\ideal{R})_k$ is the completion of
   $(L/\ideal{r})_k$, where the bracket is extended by
  continuity. In particular, if $L/\ideal{r}$ is a
  free $\integers$-module, 
  then  $\graded{F}/\ideal{R}$
  is a free $\integers_p$-module.
\end{lemma}
\begin{proof}
  The closure of $\ideal{r}$ obviously is contained in
  $\ideal{R}$. Because the bracket is continuous,
  $\graded{F}\cdot\ideal{r}=
  \overline{L}\cdot\ideal{r}\subset\overline{L\ideal{r}}\subset\overline{\ideal{r}}$.

  We have the exact sequence $0\to \ideal{r}_k\to{L}_k\to
  (L/\ideal{r})_k\to 0$. E.g.~by \cite[Lemma
  3.6]{Linnell-Schick(2000b)} the pro-$p$ completion gives rise to
  a new exact sequence
  \begin{equation*}
    0\to \ideal{R}_k\to \graded{F}_k\to (\graded{F}/\ideal{R})_k
    \to 0.
  \end{equation*}
  Here we made use of the fact that we know the completion of $L_k$
  and the closure of $\ideal{r}_k$. The assertion of the cited lemma
  of course is that the quotient $(\graded{F}/\ideal{R})_k$ is the
  completion of $(L/\ideal{r})_k$. The bracket is known to be
  continuous, therefore it extends by continuity.
\end{proof}

Let $U$ be the enveloping algebra of $L/\ideal{r}$ and $V$ be the
enveloping algebra of $\graded{F}/\ideal{R}$. Then
$\ideal{r}/[\ideal{r},\ideal{r}]$ is a $U$-module via the
adjoint action, and $\ideal{R}/\overline{[\ideal{R},\ideal{R}]}$
is a $V$-module via the adjoint action (the closure again is taken
degree-wise).

Assume that $L/\ideal{r}$ is a free $\integers$-module. Remember that the
enveloping algebra $U$ is a quotient of the tensor algebra
$T(L/\ideal{r})$ \cite[Chapter XIII.1]{Cartan-Eilenberg(1956)}. The
algebra $T(L/\ideal{r})$  has a grading induced from the grading of
$L/\ideal{r}$, where the degree of a product $x_1\tensor
x_2\tensor\dots\tensor
x_\nu$ of homogeneous elements of $L/\ideal{r}$ is the sum of the
degrees of the tensor factors $x_1,\dots,x_\nu$.  Since the relations
$x\tensor y-y\tensor x=[x,y]$ ``preserve'' the grading, it
passes to the quotient algebra $U$.

Assume that $(L/\ideal{r})_k$ is freely generated as $\integers$-module
by the finitely many generators $x_{k,1},\dots,x_{k,n_k}$. We have
seen before that $(\graded{F}/\ideal{R})_k$ is the pro-$p$
completion of $(L/\ideal{r})_k$. Therefore it is a free
$\integers_p$-module generated by $x_{k,1},\dots,x_{k,n_k}$. Observe
that the grading of $\graded{F}/\ideal{R}$ gives rise to a grading
of the enveloping algebra $V$ exactly as in the discrete case. 

In this situation, by the Poincar{\'e}-Witt theorem \cite[Theorem 3.1 in Chapter
XIII]{Cartan-Eilenberg(1956)} $U_k$ is a finitely generated free
$\integers$-module on certain products of the $x_{k,l}$, and $V_k$ is
a free $\integers_p$-module with the same basis. In particular, $V$
degree-wise is the pro-$p$ completion of $U$. Moreover,
the product $U_k\times U_l\to U_{k+l}$ extends by continuity to
$V_k\times V_l\to U_{k+l}$.

\begin{lemma}\label{lem:commutator_completion}
  In the above situation,
  $(\ideal{R}/\overline{[\ideal{R},\ideal{R}]})_k$ is for
  every $k\in\naturals$ the
  pro-$p$ completion of
  $(\ideal{r}/[\ideal{r},\ideal{r}])_k$.

  The $U$-module multiplication
  \begin{equation*}
U_k\times (\ideal{r}/[\ideal{r},\ideal{r}])_l\to (\ideal{r}/[\ideal{r},\ideal{r}])_{k+l}
\end{equation*}
extends for every $k$ and $l$ by continuity to the $V$-module multiplication
\begin{equation*}
V_k\times (\ideal{R}/\overline{[\ideal{R},\ideal{R}]})_l\to
(\ideal{R}/\overline{[\ideal{R},\ideal{R}]})_{k+l}.
\end{equation*}
\end{lemma}
\begin{proof}
  We checked that $\ideal{R}$ is (degree-wise) the closure of
  $\ideal{r}$ in $\graded{F}$. Continuity of the bracket implies
  (similar to the reasoning above) that the (degree-wise closed)
  commutator ideal $\overline{[\ideal{R},\ideal{R}]}$ is the
  (degree-wise) closure of $[\ideal{r},\ideal{r}]$. Because of
  Lemma \ref{lem:unique_completion} $\ideal{R}_k$ really is the
  pro-$p$ completion of $\ideal{r}$ for every $k$. Similar to the proof
  of Lemma \ref{lem:ideal_closure}, completing the exact sequence
  \begin{equation*}
    0\to [\ideal{r},\ideal{r}]_k\to \ideal{r}_k\to
    (\ideal{r}/[\ideal{r},\ideal{r}])_k \to 0
  \end{equation*}
  gives the exact sequence
  \begin{equation*}
    0\to \overline{[\ideal{R},\ideal{R}]}_k\to \ideal{R}_k\to
    (\ideal{R}/\overline{[\ideal{R},\ideal{R}]})_k \to 0
  \end{equation*}
  which shows that the quotient
  $\ideal{R}/\overline{[\ideal{R},\ideal{R}]}$ indeed is the
  completion of $\ideal{r}/[\ideal{r},\ideal{r}]$.

  The last statement follows immediately from continuity of the
  $V$-module multiplication map
  \begin{equation*}
V_k\times (\ideal{R}/\overline{[\ideal{R},\ideal{R}]})_l\to
(\ideal{R}/\overline{[\ideal{R},\ideal{R}]})_{k+l}.
\end{equation*}
  Continuity of this map follows from the fact that the bracket in
  $\graded{F}$ which induces this map is continuous.
\end{proof}

\begin{theorem}\label{theo:main_cond}
  Let $w_1,\dots,w_r\in F$ and $W_1,\dots,W_r\in L$ be given as above,
  with $W_1,\dots,W_r$ generating the  ideals $\ideal{r}$ of $L$ and
  $\ideal{R}$ of
  $\graded{F}$. Assume that $L/\ideal{r}$ is a free $\integers$-module
  with enveloping algebra $U$
  and $(L/\ideal{r})_k$ is finitely generated as a $\integers$-module
  for every $k$. Assume further that $\ideal{r}/[\ideal{r},\ideal{r}]$ is a
  free $U$-module with basis the images of $W_1,\dots,W_r$. in
  $\ideal{r}/[\ideal{r},\ideal{r}]$.

  Let $R$ be the closed normal subgroup of $F$ generated by
  $w_1,\dots,w_r$ and set $G=F/R$.
  Then $R/\overline{[R,R]}$ is a free $\integers_p[[G]]$-module with basis
  $w_1,\dots,w_r$.
\end{theorem}
\begin{proof}
  By Lemma \ref{lem:ideal_closure} the assumptions imply that
  $\graded{F}/\ideal{R}$ is a free $\integers_p$-module. Moreover,
  if $V$ denotes the enveloping algebra of $\graded{F}/\ideal{R}$
  then $V$ is degree-wise the pro-$p$ completion of $U$. By Lemma
  \ref{lem:commutator_completion}
  $\ideal{R}/\overline{[\ideal{R},\ideal{R}]}$ is degree-wise
  the completion of $\ideal{r}/[\ideal{r},\ideal{r}]$, which
  is by assumption a free $U$-module freely generated by the images of
  $W_1,\dots,W_r$. Since the completion of $U$ is $V$, this implies that 
  $\ideal{R}/\overline{[\ideal{R},\ideal{R}]}$ is a free
    $V$-module with basis the images of $W_1,\dots,W_r$.

    It follows that all the conditions of the pro-$p$ version of
    \cite[Theorem 1 and Theorem 2]{Labute(1985)} are fulfilled. As
    remarked at the bottom of page~52 of \cite{Labute(1985)}, the
    pro-$p$ version is obtained from the discrete version by
    everywhere replacing the ring of integers with $\integers_p$, the
    group ring $\integers[G]$ with the completed group algebra
    $\integers_p[[G]]$ and ideals and subgroups with closed ideals and
    closed subgroups. With these changes, the proofs of Theorem~1 and
    Theorem~2 of \cite{Labute(1985)} almost literally remain the same
    to give the pro-$p$ versions of these theorems (the same is true
    for the pro-$p$ versions of statements and proofs in
    \cite{Labute(1970)} which are used in \cite{Labute(1985)}). In
    particular,
    $V=\graded{\integers_p[[G]]}$ where the filtration of
    $\integers_p[[G]]$ is induced by the augmentation ideal. Moreover,
    in the course of the proof of \cite[Theorem 1]{Labute(1985)} it is
    established that the assumptions imply that
    $\mathfrak{R}/\overline{[\mathfrak{R},\mathfrak{R}]}$ is
    isomorphic as $V$-module to $\graded{R/\overline{[R,R]}}$. The
    latter therefore is a free $\graded{\integers_p[[G]]}$-module
    freely generated by the images $W_1,\dots,W_r$ of
    $w_1,\dots,w_r$. More details of the extension of the proof of
    \cite[Theorem 1]{Labute(1985)} are given in \cite[after
    Theorem 6.14]{kuempel07:_towar_atiyah}. By Lemma
    \ref{lem:graded_to_ungraded_free}
    $R/\overline{[R,R]}$ is a free $\integers_p[[G]]$-module with
    basis $w_1,\dots,w_r$.
\end{proof}

We are now in the situation to prove Theorem \ref{theorem:main}.
 Indeed, Theorem \ref{theorem:main}
follows from Proposition \ref{prop:linkpresentation} if we show that
 $R/\overline{[R,R]}$ is a free $\integers_p[[\hat G]]$-module freely
 generated by $d-1$ elements. Remember from the proof of Proposition
 \ref{prop:linkpresentation} that $R$ is the closed subgroup of
 $F:=\hat\Gamma_d$ ($\Gamma_d$ a free  group on $d$
 generators and $\hat\Gamma_d$ its pro-$p$ completion)
 that is generated by elements $w_1,\dots,w_{d-1}$ of $\hat F$ which
 are congruent to an element of $F$ modulo $\gamma_n(F)$ for every
 $n$. In particular, it follows that the images $W_1,\dots, W_{d-1}$
 of $w_1,\dots,w_{d-1}$ in $\graded{F}$ belong to the subset
 $\graded{\Gamma_d}=L$. Moreover, $\hat G=F/R$.

 Let $\mathfrak{r}$ be the ideal of $L$ generated by
 $W_1,\dots,W_{d-1}$ and let $U$ be the enveloping algebra of
 $L/\mathfrak{r}$. By \cite[Theorem 2]{Labute(1990)}, $L/\mathfrak{r}$
 is a free 
 $\integers$-module which, by \cite[Theorem 1]{Labute(1990)}, is
 finitely generated in each degree, and
 $\mathfrak{r}/[\mathfrak{r},\mathfrak{r}]$ is a free $U$-module
 freely generated by the images of $W_1,\dots,W_{d-1}$. Then Theorem
 \ref{theo:main_cond} implies that $R/\overline{[R,R]}$ indeed is a
 free $\integers_p[[\hat G]]$-module and Theorem \ref{theorem:main} follows.

\bibliographystyle{amsplain}
\bibliography{link_Galoiscohomology}

\end{document}